\documentclass[reqno]{amsart}
\usepackage{amsmath,amsthm,amssymb,amscd}
\usepackage{hyperref}

\newcommand{\F}{\mathbb F}

\newcommand{\Z}{\mathbb Z}
\newcommand{\Q}{\mathbb Q}

\newcommand{\eps}{\varepsilon}

\newcommand{\cD}{\mathcal D}
\newcommand{\cO}{\mathcal O}

\newcommand{\fp}{{\mathfrak p}}

\newcommand{\Qt}{{\mathbb Q}^\times}
\newcommand{\Qts}{\Q^{\times\,2}}
\newcommand{\im}{{\operatorname{im}}\,}
\newcommand{\Spl}{{\operatorname{Spl}}}
\newcommand{\Gal}{{\operatorname{Gal}}}

\newcommand{\disc}{{\operatorname{disc}}\,}
\newcommand{\rank}{{\operatorname{rank}}\,}

\newcommand{\Cl}{{\operatorname{Cl}}}

\newfont{\cyr}{wncyb10}

\newcounter{lemmacount}[section]

\newtheorem{thm}[lemmacount]{Theorem}
\newtheorem{prop}[lemmacount]{Proposition}
\newtheorem{lem}[lemmacount]{Lemma}
\newtheorem{cor}[lemmacount]{Corollary}

\title[Higher Descent on Pell Conics]{Higher Descent on Pell Conics. \\
      I. From Legendre to Selmer}
\author{Franz Lemmermeyer}
\address{Department of Mathematics,
  Bilkent University,
  06800 Bilkent, Ankara, Turkey}
\email{franz@fen.bilkent.edu.tr}

\begin{document}
\maketitle

\section*{Introduction}

The theory of Pell's equation has a long history, as can be seen
from the huge amount of references collected in Dickson \cite{Dick},
from the two books on its history by Konen \cite{Kon} and Whitford
\cite{Whit}, or from the books by Walfisz \cite{Wal}, Faisant \cite{Fai}, 
and Barbeau \cite{Bar}. For the better part of the last few centuries, 
the continued fractions method was the undisputed method for 
solving a given Pell equation, and only recently faster methods have 
been developed (see the surveys by Lenstra \cite{Len} and H.C.~Williams 
\cite{Wil}).

This is the first in a series of articles which has the goal of
developing a theory of the Pell equation that is as close as possible
to the theory of elliptic curves: we will discuss $2$-descent on Pell 
conics, introduce Selmer and Tate-Shafarevich groups, and formulate
an analogue of the Birch--Swinnerton-Dyer conjecture. 

In this article, we will review the history of results that are
related to this new interpretation. We will shortly discuss the
construction of explicit units in quadratic number fields, and then
deal with Legendre's equations and the results of R\'edei, Reichardt
and Scholz on the solvability of the negative Pell equation.

The second article \cite{L2} is devoted to references to a 
``second $2$-descent'' in the mathematical literature from Euler 
to our times, and in \cite{L3} we will discuss the first 
$2$-descent and the associated Selmer and Tate-Shafarevich
groups from the modern point of view.

\section{Explicit Units}

Since finding families of explicit units in number fields
is only indirectly related to our topic, we will be rather 
brief here. The most famous families of explicitly given 
units live in fields of Richaud-Degert type: if $d = b^2 + m$ 
and $b \mid 4m$, then $\alpha = b+\sqrt{d}$ has norm $-m$, 
hence is a unit if $m \in \{\pm 1, \pm 4\}$; if $\alpha$ 
is not a unit, then $\frac1m \alpha^2$ is. 

Brahmagupta observed in 628 AD that if $a^2 - nb^2 = k$, then 
$x = \frac{2ab}k$, $y = \frac{a^2 + nb^2}k$ satisfy
the Pell equation $x^2 - ny^2 = 1$. If $k = \pm 1, \pm 2$,
these solutions are integral; if $k = \pm 4$, Brahmagupta
showed how to construct an integral solution. Note that
$x$ and $y$ are necessarily integral if $k$ divides the
squarefree integer $n$.
 
According to C.~Henry (see \cite{Male} and Dickson 
\cite[v. II, p. 353]{Dick}), Malebranche (1638--1715) 
claimed that the Pell equation $Ax^2 + 1 = y^2$ can be 
solved easily if $A = b^2 + m$ with $m|2b$.

Euler mentioned in \cite[p. 99]{Eul4} and later again in 
\cite{Euler} that if $d = b^2c^2 \pm 2b$ or $d = b^2c^2 \pm b$, 
then the Pell equation $x^2 - dy^2 = 1$ can be solved explicitly. 
This result was rediscovered e.g. by M.~Stern \cite{Stern}, 
Richaud \cite{Ri3}, Hart \cite{Hart2}, Speckman \cite{Spec}
and Degert \cite{Deg}. 
Special cases are due to Moreau \cite{Mor73}, 
de Jonqui\`eres \cite{Jon},Ricalde \cite{Ric}, 
Boutin \cite{Bout}, Malo \cite{Mal}, and von Thielmann \cite{vTh}. 
The quadratic fields $\Q(\sqrt{d}\,)$ with $d = a^2 \pm r$ 
and $r \mid 2a$ were called fields of Richaud-Degert type 
by Hasse \cite{Has}. Degert's results were generalized by
Yokoi \cite{Yok,Yok2}, Kutsuna \cite{Kut}, Katayama 
\cite{Kat1,Kat2}, Takaku \& Yoshimoto \cite{TY}, Ramasamy 
\cite{Ram}, and Mollin \cite{Mol1}.

The results about units in ``fields of Richaud-Degert type''
had actually been generalized long before Degert. Observe 
that the computation of the fundamental solution of 
$X^2 - dY^2 = 1$ for $d = m^2+1$ with the method of
continued fractions is trivial because the development
of $\sqrt{d}$ has period $1$. Similarly, the development
of $\sqrt{d}$ for $d$ of Richaud-Degert type have small
periods; here are some examples:
\begin{itemize}
\item $d = k^2+k$: $\sqrt{d} = [k,\overline{2,2k}]$;
\item $d = k^2+2k$: $\sqrt{d} = [k,\overline{1,2k}]$;
\item $d = a^2k^2+a$: $\sqrt{d} = [ak,\overline{2k,2ak}]$;
\item $d = a^2k^2+2a$: $\sqrt{d} = [ak,\overline{k,2ak}]$.
\end{itemize}
Actually, already Euler \cite{Eul4} gave the continued fraction
expansions for $d = n^2+1$, $n^2+2$, $n^2+n$, $9n^2+3$ and a few
other values of $d$.

Perron \cite{Perr} gave examples of polynomials $f(x)$ for 
which the continued fraction expansion of $\sqrt{f(x)}$ can be 
given explicitly; the period lengths of his examples were 
$\le 6$. More examples were given by Yamamoto \cite{Yam} 
and Bernstein \cite{Ber1,Ber2} (who produced units whose 
continued fraction expansions have arbitrarily long period), 
as well as by Azuhata \cite{Azh1,Azh2}, Tomita \cite{Tom1,Tom2}, 
Levesque \& Rhin \cite{LR}, Levesque \cite{Lev}, Madden \cite{Mad}, 
Mollin \cite{Mol}, H.C.~Williams \cite{WilS},
 van der Poorten \& H.C.~Williams \cite{PW}, 
J. McLaughlin \cite{McL1,McL2}, and probably many others.

Nathanson \cite{Nat} proved that $X^2 - DY^2 = 1$, where
$D= x^2+d$, has nontrivial solutions $X, Y \in \Z[x]$ if and
only if $d = \pm 1, \pm 2$. This result was generalized by 
Hazama \cite{Haz} and by Webb \& Yokota \cite{WY}. For
connections with elliptic curves, see 
Berry \cite{Ber} and Avanzi \& Zannier \cite{AZ}.

\section{Legendre}

\subsection{Legendre's Th\'eorie des Nombres}
\S\ VII of Legendre's book \cite{LegTN} on number theory had the
title
\begin{quote}
Th\'eor\`emes sur la possibilit\'e des \'equations de la forme
$Mx^2 - Ny^2 = \pm 1$ ou $\pm 2$.\footnote{Theorems on the
solvability of the equations of the form $Mx^2 - Ny^2 = \pm 1$ 
or $\pm 2$.}
\end{quote}

Legendre starts his investigation by assuming that $A$ is prime, 
and that $p$ and $q$ are the smallest positive solutions of the 
equation 
\begin{equation}\label{ELe}
p^2 - Aq^2 = 1.
\end{equation} 
Writing this equation as $(p-1)(p+1) = Aq^2$ he deduces that 
$q = fgh$ with $f \in \{1, 2\}$ and 
$$ \left. \begin{array}{rcl}
     p+1 & = & fg^2A \\  p-1 & = & fh^2 \end{array} \right\}
   \quad \text{or} \quad
   \left. \begin{array}{rcl}
     p+1 & = & fg^2 \\  p-1 & = & fh^2A \end{array} \right\} $$

Subtracting these equations from each other he gets equations of
the form $\pm \frac2f = x^2 - Ay^2$. The case $f = \pm 1$ leads
to contradictions modulo $4$, the case $f = 2$ contradicts the
minimality of $p$ and $q$, thus we must have $f = -2$ and therefore
$-1 = x^2 - Ay^2$ (see \cite[p. 55]{LegTN}).
 
Similar arguments easily yield

\begin{prop}\label{PL1}
Let $p$ be a prime.
\begin{enumerate}
\item If $p \equiv 1 \bmod 4$, then $X^2 - pY^2 = -1$ 
      has integral solutions.
\item If $p \equiv 3 \bmod 8$, then $X^2 - pY^2 = -2$ 
      has integral solutions.
\item If $p \equiv 7 \bmod 8$, then $X^2 - pY^2 = +2$ 
      has integral solutions.
\end{enumerate}
\end{prop}

Next Legendre considers composite values of $A$: if $A = MN$ 
is the product of two odd primes $\equiv 3 \bmod 4$, then he 
shows that $Mx^2 - Ny^2 = \pm 1$ is solvable;\footnote{Using 
the solvability of $Mx^2 - Ny^2 = 1$ for primes $M \equiv N 
\equiv 3 \bmod 4$, Legendre later proved a special case of the 
quadratic reciprocity law, namely that $(\frac{p}q) = - (\frac{q}p)$ 
for such primes.} 
if $M \equiv N \equiv 1 \bmod 4$, however, none of the equations 
$x^2 - MNy^2 = -1$ and $Mx^2 - Ny^2 = \pm 1$ can be excluded. He 
also states that for given $A$ at most one of these equations can
have a solution, but the argument he offers is not conclusive.
 
He then treats the general case and comes to the conclusion
\begin{quote}
\'Etant donn\'e un nombre quelquonque non quarr\'e $A$, il est toujours
possible de d\'ecomposer ce nombre en deux facteurs $M$ et $N$ tels que
l'une des deux \'equations $Mx^2 - Ny^2 = \pm 1$, $Mx^2 - Ny^2 = \pm 2$
soit satisfaite, en prenant convenablement le signe du second 
membre.\footnote{Given an arbitrary nonsquare number $A$ it is always
possible to decompose this number into two factors $M$ and $N$ such that
one of the two equations $Mx^2 - Ny^2 = \pm 1$, $Mx^2 - Ny^2 = \pm 2$
is satisfied, the sign of the second number being conveniently taken.} 
\end{quote} 

Again, Legendre adds the remark that among these equations there is
exactly one which is solvable; he gives a different argument this time,
which again is not conclusive since he is only working with the minimal
solution of the Pell equation.

\subsection{Dirichlet's Exposition}
Dirichlet \cite{Dir} gave an exposition of Legendre's technique,
which we will quote now; his \S\ 1 begins 
\begin{quote}
Wir beginnen mit einer kurzen Darstellung der Legendre'schen Methode. 
Es bezeichne $A$ eine gegebene positive Zahl ohne qua\-dra\-ti\-schen 
Factor, d.h. deren Primfactoren alle von einander verschieden sind,
und es seien $p$ und $q$ die kleinsten Werthe $(p = 1$ und $q = 0$
ausgenommen), welche der bekanntlich immer l\"osbaren Gleichung:
\begin{equation}\label{ED1} \tag{1} p^2 - Aq^2 = 1 \end{equation}
Bringt man dieselbe in die Form $(p+1)(p-1) = Aq^2$, und bemerkt
man, dass $p+1$ und $p-1$ relative Primzahlen sind oder bloss
den gemeinsamen Factor $2$ haben, je nachdem $p$ gerade oder
ungerade ist, so sieht man, dass die Gleichung (\ref{ED1}) im
ersten Falle die folgenden nach sich zieht:
$$ p+1 = Mr^2, \quad p-1 = Ns^2, \quad A = MN, \ q = rs, $$
und ebenso im zweiten:
$$ p+1 = 2Mr^2, \quad p-1 = 2Ns^2, \quad A = MN, \ q = 2rs, $$
wo $M, N$ und mithin $r, s$ durch $p$ v\"ollig bestimmt sind. Es
sind n\"amlich $M$, $N$ im ersten Fall respective die gr\"ossten
ge\-mein\-schaft\-li\-chen Theiler von $A$, $p+1$ und $A$, $p-1$, im 
andern dagegen von $A$, $\frac{p+1}2$ und $A$, $\frac{p-1}2$.
Aus diesen Gleichung folgt 
\begin{equation}\label{ED2}
\tag{2} Mr^2- Ns^2 = 2, \quad Mr^2- Ns^2 = 1.
\end{equation}

Hat man die Gleichung (\ref{ED1}) nicht wirklich aufgel\"ost, und
ist also $p$ nicht bekannt, so weiss man bloss, dass eine dieser
Gleichungen stattfinden muss, und da unter dieser Voraussetzung 
$M$ und $N$ nicht einzeln gegeben sind, so enth\"alt jede der 
Gleichungen (\ref{ED2}) mehrere besondere Gleichungen, die man
erh\"alt, indem man successive f\"ur $M$ alle Factoren von $A$
($1$ und $A$ mit eingeschlossen) nimmt und $N = \frac{A}{M}$ setzt.
\footnote{We start with a short presentation of Legendre's method.
Let $A$ denote a given positive number without quadratic factors,
i.e., whose prime factors are all pairwise distinct, and let
$p$ and $q$ be the smallest values (except $p = 1$ and $q = 0$)
satisfying the equation (\ref{ED1}), which is known to be
always solvable. Writing (\ref{ED1}) in the form 
$(p+1)(p-1) = Aq^2$ and observing that $p+1$ and $p-1$ are 
relatively prime or have the common factor $2$ according as
$p$ is even or odd, it can be seen that the equation (\ref{ED1}) 
in the first case implies the following:
$$ p+1 = Mr^2, \quad p-1 = Ns^2, \quad A = MN, \ q = rs, $$
and similarly in the second case:
$$ p+1 = 2Mr^2, \quad p-1 = 2Ns^2, \quad A = MN, \ q = 2rs, $$
where $M, N$ and therefore $r, s$ are completely determined by $p$.
In fact, $M$, $N$ are, in the first case, the greatest common divisors
of $A$, $p+1$ and $A$, $p-1$, in the second case of $A$, $\frac{p+1}2$ 
and $A$, $\frac{p-1}2$. From these equations we deduce (\ref{ED2}).

If we do not have a solution of (\ref{ED1}), thus if $p$ is not known,
then we only know that one of these equations must have a solution,
and since under this assumption $M$ and $N$ are not given explicitly,
each of the equations (\ref{ED2}) contains several equations which
we get by letting $M$ run through all factors of $A$ (including $1$ 
and $A$) and putting $N = \frac{A}{M}$.}
\end{quote}
Thus we have

\begin{prop}
The solvability of the Pell equation $X^2 - AY^2 = 1$ implies the 
solvability of one of the equations $Mr^2 - Ns^2 = 1$ or 
$Mr^2 - Ns^2 = 2$, where $MN = A$.
\end{prop}

This set of auxiliary equations (\ref{ED2}) was derived using 
continued fractions by Arndt \cite{Arn1,Arn2}, and later also
by Richaud \cite{Ri2,Ri3} and Roberts \cite{Rob}.
Catalan \cite{Cat} also studied the equations (\ref{ED2}). 

Legendre's method, Dirichlet writes, consists
in showing that all but one of these equations are
unsolvable, thus demonstrating that the remaining 
equation must have integral solutions. Dirichlet then
shows how to apply this technique by proving Legendre's
Proposition \ref{PL1} as well as the following result:

\begin{prop}\label{PL2}
Let $p$ denote a prime.
\begin{enumerate}
      \item If $p \equiv 3 \bmod 8$, then $pr^2 - 2s^2 = 1$ 
            has integral solutions. 
      \item If $p \equiv 5 \bmod 8$, then $2pr^2 - s^2 = 1$ 
            has integral solutions.
      \item If $p \equiv 7 \bmod 8$, then $2r^2 - ps^2 = 1$ 
            has integral solutions. 
\end{enumerate}
If $p \equiv 1 \bmod 8$, either of the three equations
may be solvable.
\end{prop}

\subsection{Richaud}
These results were extended to $d$ having many prime factors
by Richaud. In \cite{Ri1}, he states without proof some results 
of which the following are special cases (Richaud considered 
values of $d$ that were not necessarily squarefree):

\begin{prop}
\begin{enumerate}
\item If $p$ and $q$ are primes congruent to $5 \bmod 8$, then the
      equations $X^2 - 2pY^2 = -1$ and $X^2 - 2pqY^2 = -1$ are
      solvable in integers.
\item If $p, q$ and $r$ are primes congruent to $5 \bmod 8$, and
      if $(2p/q) = (2p/r) = -1$, then the equation 
      $X^2 - 2pqrY^2 = -1$ is solvable in integers.
\item If $p \equiv 5 \bmod 8$ and $a \equiv b \equiv 1 \bmod 8$
      are primes such that $(2p/a) = (2p/b) = -1$, then the
      equations $X^2 - 2apY^2 = -1$ and $X^2 - 2abpY^2 = -1$ are
      solvable in integers.
\item If $p \equiv q \equiv 5 \bmod 8$ and $a \equiv 1 \bmod 8$
      are primes such that $(pq/a) = -1$, then the equation 
      $X^2 - 2apqY^2 = -1$ is solvable in integers.
\end{enumerate}
\end{prop}

Here are some results from Richaud \cite{Ri2,Ri3,Ri4}:

\begin{prop}
Let $p_i$ denote primes. 
\begin{enumerate}
\item If $d = p_1p_2p_3p_4$, $p_i \equiv 5 \bmod 8$, if
      $X^2 - p_1p_2p_3 Y^2 = -1$, and if 
      $(p_4/p_1) = (p_4/p_2) = (p_4/p_3) = 1$, then $X^2 - dY^2 = -1$ 
      is solvable. This generalizes to arbitrary numbers of primes.
\item If $d = p_1p_2p_3p_4$, $p_i \equiv 5 \bmod 8$, if
      $X^2 - p_1p_2p_3 Y^2 = -1$, and if 
      $(p_1/p_3) = (p_2/p_3) = -1$, $(p_1/p_4) = (p_2/p_4) = (p_3/p_4) = -1$, 
      then $X^2 - dY^2 = -1$ is solvable. This generalizes to arbitrary 
      numbers of primes.
\item If $d = p_1p_2p_3p_4$, where 
      $p_1 \equiv \ldots \equiv p_4 \equiv 5 \bmod 8$, then 
      $X^2 - dY^2 = -1$ is solvable if 
      $(p_1p_2/p_3) = (p_1p_2/p_4) = (p_3p_4/p_1) = (p_3p_4/p_2)$.
\end{enumerate}
\end{prop}

The following result is due to Tano \cite{Tano}; the 
special case $n = 3$ was already proved by Dirichlet.

\begin{prop}\label{PTa}
Let $p_1$, \ldots, $p_n$ denote primes $p_i \equiv 1 \bmod 4$,
where $n \ge 3$ is an odd number, and put $d = p_1\cdots p_n$. 
Assume that $(p_i/p_j) = +1$ for at most one pair $(i,j)$. 
Then $X^2 - dY^2 = -1$ has an integral solution.
\end{prop}

Proposition \ref{PTa} was generalized by Trotter \cite{Tro},
who was motivated by the results of Pumpl\"un \cite{Pum}:

\begin{prop}
Let $p_1$, \ldots, $p_n$, where $n$ is an odd integer, denote 
primes $p_i \equiv 1 \bmod 4$. If there is no triple $i, j, k$ 
such that $(p_i/p_j) = (p_j/p_k) = +1$, then $X^2 - dY^2 = -1$ 
has an integral solution.
\end{prop}

Newman \cite{New} rediscovered a weaker form of Proposition \ref{PTa}: 
he assumed that  $(p_i/p_j) = -1$ for all $i \ne j$. 

There are a lot more results of this kind concerning the \
solvability of the negative Pell equation $X^2 - dY^2 = -1$ 
in terms of Legendre symbols than we can (or may want) to 
mention here. We will see in \cite{L3} that these results 
can be interpreted as computations of Selmer groups for
certain types of discriminants.

In fact, as Dickson \cite[\S 25]{DiSN} observed, Legendre's set 
of equations $Mr^2 - Ns^2 = 1, 2$ admit a group structure; in
modern language, it is induced by identifying $Mr^2 - Ns^2 = 1$ 
and $Mr^2 - Ns^2 = 2$ with the elements $M\Qts$ and $2M\Qts$ of 
the multiplicative group $\Qt/\Qts$. This group will be called 
the $2$-Selmer group of the corresponding Pell equation (see \cite{L3}). 

\section{Dirichlet}

After having explained Legendre's technique, Dirichlet \cite{Dir}
refined this method by invoking the quadratic reciprocity law. 
His first result going beyond Legendre is the following (for
information on quartic residue symbols see \cite{LRL}):

\begin{prop}\label{PDi1}
For primes $p \equiv 9 \bmod 16$ such that $(2/p)_4 = -1$,
the equation $2pr^2 - s^2 = 1$ has an integral solution.
\end{prop}

\begin{proof}
Consider the equation $2r^2 - ps^2 = 1$. We see that $s$ is odd, 
and that $(2/s) = +1$. Thus $s \equiv \pm 1 \bmod 8$, therefore
$s^2 \equiv1 \bmod 16$, and then $2r^2 \equiv 2 \bmod 16$
and $p \equiv 9 \bmod 16$ yield a contradiction.

Next consider $pr^2 - 2s^2 = 1$. Here $1 = (2/p)_4(s/p)$,
and $(s/p) = (p/s) = 1$. Thus solvability implies $(2/p)_4 = 1$,
which is a contradiction.
\end{proof}

Dirichlet next considers other cases where $d$ has two 
or three prime factors; we are content with mentioning 

\begin{prop}\label{PDi2}
Let $p \equiv q \equiv 1 \bmod 4$ be distinct primes. 
\begin{enumerate}
\item If $(p/q) = -1$, then $s^2 - pqr^2 = -1$ has an integral solution.
\item If $(p/q)_4 = (q/p)_4 = -1$, then $s^2 - pqr^2 = -1$ has an 
      integral solution.
\end{enumerate}
\end{prop}

Using similar as well as some other techniques, Dirichlet's 
results such as Propositions \ref{PDi1} and \ref{PDi2} were 
generalized by Tano \cite{Tano}.

Dirichlet was the first who observed that there is essentially
only one among the equations derived by Legendre which has an 
integral solution:

\begin{thm}\label{TDi0}
Let $A$ be a positive squarefree integer. Then there is exactly 
one pair of positive integers $(M,N) \ne (1,A)$ with $MN = d$ 
such that $Mr^2 - Ns^2 = 1$ or $Mr^2 - Ns^2 = 2$ is solvable.
\end{thm}

\begin{proof}[Proof (Dirichlet).]
We know that all solutions of the Pell equation $P^2 - AQ^2 = 1$ 
are given by $P + Q\sqrt{A} = \pm (p+q\sqrt{A})^m$, where $m \in \Z$
and where $p+q\sqrt{A}$ is the fundamental solution. Up to sign
we thus have 
$ P = \frac12((p+q\sqrt{A})^m + (p-q\sqrt{A})^m). $
This shows that $P \equiv p^m \bmod A$.

Assume first that $m$ is odd. Then $P \equiv p \bmod 2$. If $p$ 
is even, then we have 
$$ p \equiv -1 \bmod M, \qquad p \equiv 1 \bmod N,$$ 
which implies that 
$$ P \equiv -1 \bmod M, \qquad P \equiv 1 \bmod N.$$
Thus $P+Q\sqrt{A}$ leads to the very same equation
$Mr^2 - Ns^2 = 2$ as the fundamental solution $p+q\sqrt{A}$.
If $p$ is odd,  then we find similarly that
$$ p \equiv -1 \bmod 2M, \qquad p \equiv 1 \bmod 2N,$$ 
and again $P+Q\sqrt{A}$ leads to the same equation
$Mr^2 - Ns^2 = 1$ as $p+q\sqrt{A}$.

Finally, if $m$ is even, then it is easy to see that 
$P+Q\sqrt{A}$ leads to $r^2 - As^2 = 1$. This shows that
there are exactly two equations with integral solutions.
\end{proof}

Proofs of results equivalent to Theorem \ref{TDi0}, based on
the theory of continued fractions, were given by Petr 
\cite{Pet1,Pet2}, and Halter-Koch \cite{HK};
different proofs are due to Nagell \cite{Nag3}, 
Kaplan \cite{Kap}, and Mitkin \cite{Mit}; Trotter \cite{Tro} 
proved the special case where the fundamental unit has negative 
norm. Pall \cite{Pal} showed that Theorem \ref{TDi0} follows 
easily from a result that Gauss proved in his Disquisitiones 
Arithmeticae. See also Walsh \cite{Walth,Walsh}.

In the ideal-theoretic interpretation, Theorem \ref{TDi0} says
that if $K = \Q(\sqrt{d}\,)$ has a fundamental unit of norm $+1$,
then there is exactly one nontrivial relation in the usual class 
group among the ambiguous ideals, the trivial relation coming from
the factorization of the principal ideal $(\sqrt{d}\,)$. From this
point of view, Theorem \ref{TDi0} is an important step in the proof
of the ambiguous class number formula for quadratic number fields.

\section{R\'edei, Reichardt, Scholz}

In 1932, R\'edei started studying the $2$-class group of the 
quadratic number field $k = \Q(\sqrt{d}\,)$, and applied it in 
\cite{RedPl} to problems concerning the solvability of the negative 
Pell equation $X^2 - dY^2 = -4$. In the following, let $e_2$ and 
$e_4$ denote the $2$-rank and the $4$-rank of the class group 
$C = \Cl_2^+(k)$ of $k$ in the strict sense, that is, put 
$e_2 = \dim_{\F_2} C/C^2$ and $e_4 = \dim_{\F_2} C^2/C^4$. 
A factorization of the discriminant $d = \disc k$ into 
discriminants $d = \Delta_1 \Delta_2$ is called a splitting
of the second kind (or $C_4$-decomposition) if 
$(\Delta_1/p_2) = (\Delta_1/p_2) = +1$
for all primes $p_i \mid \Delta_i$. The main result of 
\cite{Red32} (see also \cite{RR33}) is

\begin{prop}
The $4$-rank $e_4$ equals the number of independent splittings
of the second kind of $d$.
\end{prop}

This result of R\'edei turned out to be very attractive; new variants
of proofs were given by Iyanaga \cite{Iya}, Bloom \cite{Blo}, 
Carroll \cite{Car}, and Kisilevsky \cite{Kis}. Inaba \cite{Ina}, 
Fr\"ohlich \cite{Froa} and G.~Gras \cite{Gr73} investigated the 
$\ell$-class group of cyclic extensions of prime degree $\ell$; 
see Stevenhagen \cite{Ste} for a modern exposition. For
generalization's of R\'edei's technique to quadratic extensions 
of arbitrary number fields see G.~Gras \cite{Gr92}. 
Morton \cite{Mor0} and Lagarias \cite{Laga} gave modern accounts
of R\'edei's method for computing the $2$-part of the class
groups of quadratic number fields.

Damey \& Payan \cite{DP} proved

\begin{prop}
Let $k^+ = \Q(\sqrt{m}\,)$ be a real quadratic number field, 
and put $k^- = \Q(\sqrt{-m}\,)$. Then the $4$-ranks $r_4^+(k^+)$ 
and $r_4(k^-)$ of $\Cl^+(k^+)$ (the class group of $k^+$ in the 
strict sense) and $\Cl(k^-)$ satisfy the inequalities 
$r_4^+(k^+) \le r_4(k^-) \le r_4^+(k^+) + 1. $
\end{prop}

Other proofs were given by G.~Gras \cite{Gr73}, Halter-Koch 
\cite{HK84}, Uehara \cite{Ueh89} and Sueyoshi \cite{Sue95}; 
see also Sueyoshi \cite{Sue97,Sue00}. Bouvier \cite{Bou1,Bou2} 
proved that the $4$-ranks and $8$-ranks of $\Q(\sqrt{2},\sqrt{m}\,)$ 
and $\Q(\sqrt{2},\sqrt{-m}\,)$ differ at most by $4$. This was 
generalized considerably by Oriat \cite{Ori,Oria,Orib}; the 
following proposition is a special case of his results:

\begin{prop}
Let $k$ be a totally real number field with odd class number,
let $2^m \ge 4$ be an integer with the property that $k$ 
contains the maximal real subfield of the field of $2^m$-th 
roots of unity, and let $d \in k^\times$ be a nonsquare. 
Then the $2^m$-ranks $r_m(K)$ and $r_m(K')$ of the class 
groups in the strict sense of $K = k(\sqrt{d}\,)$ and 
$K' = k(\sqrt{-d}\,)$ satisfy
$$ r_m(K) - r_m(K') \le R^- - r,$$ 
where $R^-$ and $r$ denote the unit rank of $K^-$ and $k$, 
respectively.
\end{prop}

R\'edei \& Reichardt \cite{RR33,Red53} and Iyanaga \cite{Iya} 
observed the following 

\begin{prop}\label{ProEl}
If $\Cl_2^+(k)$ is elementary abelian, then $N\eps = -1$.
\end{prop}

Since $e_4 = 0$ is equivalent to the rank of $R(d)$ being maximal,
this can be expressed by saying that if the R\'edei matrix
has maximal possible rank $n-1$, then $N\eps = -1$.

R\'edei introduced what is now called the R\'edei matrix of
the quadratic field with discriminant $d = d_1 \cdots d_n$.
It is defined as the $n \times n$-matrix $R(d) = (a_{ij})$
with $a_{ij} = (d_i/p_j)$ for $i \ne j$ and 
$a_{ii} = \sum_{i \ne j} a_{ij}$. R\'edei proved that the
$4$-rank of $\Cl^+(k)$ is given by
$$ e_4 = n-1- \rank R(d).$$

R\'edei (see e.g. \cite{Red37}) 
introduced a group structure on the set of splittings
of the second kind by taking the product of two such factorizations
$d = \Delta_1 \Delta_2$ and $d = \Delta_1' \Delta_2'$ to be
the factorization $d = \Delta_1'' \Delta_2''$, where
$$\Delta_1'' = \frac{\Delta_1 \Delta_1'}{\gcd(\Delta_1,\Delta_1')^2}. $$
This product is well defined because of the relations
$$ \frac{\Delta_1 \Delta_1'}{\gcd(\Delta_1,\Delta_1')^2} = 
   \frac{\Delta_2 \Delta_2'}{\gcd(\Delta_2,\Delta_2')^2}, \quad
   \frac{\Delta_1 \Delta_2'}{\gcd(\Delta_1,\Delta_2')^2} = 
   \frac{\Delta_2 \Delta_1'}{\gcd(\Delta_2,\Delta_1')^2}. $$

R\'edei \& Reichardt \cite{RR33} and Scholz \cite{Sch34} started 
applying class field theory to finding criteria for the 
solvability of the negative Pell equation $X^2 - dY^2 = -4$.
The immediate connection is provided by the following simple
observation: we have $N\eps = -1$ if and only if $\Cl(K) \simeq \Cl^+(K)$,
which in turn is equivalent to the fact that the Hilbert class fields in 
the strict (unramified outside $\infty$) and in the usual sense
(unramified everywhere) coincide. This proves

\begin{prop}
The fundamental unit of the quadratic field $K$ has negative norm
if and only if the Hilbert class field in the strict sense is
totally real.
\end{prop}

The following theorem summarizes the early results of R\'edei and
Scholz; observe that `unramified' below means `unramified outside 
$\infty$':

\begin{thm}\label{RedRe}
Let $k$ be a quadratic number field with discriminant $d$. There
is a bijection between unramified cyclic $C_4$-extensions and
$C_4$-factorizations of $d$. 

If $K/k$ is an unramified $C_4$-extension, then $K/\Q$ is normal 
with $\Gal(K/\Q) \simeq D_4$. The quartic normal extension $F/\Q$
contained in $K$ can be written in the form 
$F = \Q(\sqrt{\Delta_1},\sqrt{\Delta_2}\,)$. A careful examination 
of the decomposition and inertia groups of the ramifying primes 
shows that $(\Delta_1, \Delta_2) = 1$ and that 
$d = \Delta_1\cdot \Delta_2$ is a $C_4$-factorization.

Conversely, if $d = \Delta_1 \Delta_2$ is a $C_4$-factorization 
of $d$, then the diophantine equation 
$ X^2 - \Delta_1Y^2 = \Delta_2Z^2$ has a nontrivial solution 
$(x,y,z)$, and the extension $K = k(\sqrt{\Delta_1}, \sqrt{\mu}\,)$, where
$\mu = x+y\sqrt{\Delta_1}$, is a $C_4$-extension of $k$ unramified outside
$2\infty$. By choosing the signs of $x,y,z$ suitably one can make
$K/k$ unramified outside $\infty$.  
\end{thm}

The question of whether the cyclic quartic extension $K/k$ constructed 
in Theorem 13.1. is real or not was answered by Scholz \cite{Sch34}. 
Clearly this question is only interesting if both $\Delta_1$ and $\Delta_2$ 
are positive. Moreover, if one of them, say
$\Delta_1$, is divisible by a prime $q \equiv 3 \bmod 4$, then there
always exists a real cyclic quartic extension $K/k$: this
is so because $\alpha = x+y\sqrt{\Delta_1}$ as constructed above 
is either totally positive or totally negative (since it has 
positive norm), hence either $\alpha \gg 0$ or $-q\alpha \gg 0$,
so either $k(\sqrt{\alpha}\,)$ or $k(\sqrt{-q\alpha}\,)$ is
the desired extension. We may therefore assume that
$d$ is not divisible by a prime $q \equiv 3 \bmod 4$, i.e.
that $d$ is the sum of two squares. Then Scholz \cite{Sch34}
has shown 

\begin{prop}
Let $k$ be a real quadratic number field with discriminant $d$, and
suppose that $d$ is the sum of two squares. Assume moreover that
$d = \Delta_1\cdot \Delta_2$ is a $C_4$-factorization. Then the
cyclic quartic $C_4$-extensions $K/k$ containing 
$\Q(\sqrt{\Delta_1},\sqrt{\Delta_2}\,)$ are real if and only if 
$(\Delta_1/\Delta_2)_4 (\Delta_2/\Delta_1)_4 = +1$. Moreover, if there exists
an octic cyclic unramified extension $L/k$ containing $K$,
then $(\Delta_1/\Delta_2)_4 = (\Delta_2/\Delta_1)_4$.
\end{prop}

If $\Delta_1$ and $\Delta_2$ are prime, we can say more (\cite{Sch34}):

\begin{thm} 
Let $k = \Q(\sqrt{d}\,)$ be a real quadratic number field, and 
suppose that $d = \disc k = \Delta_1\Delta_2$ is the product of two positive 
prime discriminants $\Delta_1, \Delta_2$. Let $h(k)$, $h^+(k)$ and $\eps$
denote the class number, the class number in the strict sense, and
the fundamental unit of $\cO_k$, respectively; moreover, let $\eps_1$
and $\eps_2$ denote the fundamental units of $k_1 = \Q(\sqrt{\Delta_1}\,)$
and $k_2 = \Q(\sqrt{\Delta_1}\,)$. There are the following possibilities:
\begin{enumerate}
\item $(\Delta_1/\Delta_2) = -1$: then $h(k) = h^+(k) \equiv 2 \bmod 4$,
	and $N\eps = -1$. 
\item $(\Delta_1/\Delta_2) = +1$: then $(\eps_1/\Delta_2) = (\eps_2/\Delta_1) 
	= (\Delta_1/\Delta_2)_4  (\Delta_2/\Delta_1)_4$, and \ref{RedRe}
	shows that there is a cyclic quartic subfield
	$K$ of $k^1$ containing $k_1k_2$;
	
	\noindent i) $(\Delta_1/\Delta_2)_4 = -(\Delta_2/\Delta_1)_4$: then 
	$h^+(k) = 2\cdot h(k) \equiv 4 \bmod 8$, $N\eps = +1$,
	and $K$ is totally complex;

	\noindent ii) $(\Delta_1/\Delta_2)_4 = (\Delta_2/\Delta_1)_4 = -1$:
        then $h^+(k) = h(k) \equiv 4 \bmod 8$, $N\eps = -1$,
	and $K$ is totally real.

	\noindent iii) $(\Delta_1/\Delta_2)_4 = (\Delta_2/\Delta_1)_4 = +1$: 
        then $h^+(k) \equiv 0 \bmod 8$, and $K$ is totally real.
\end{enumerate}
\end{thm}
Here $(\Delta_1/\Delta_2)_4$ denotes the rational biquadratic residue symbol
(multiplicative in both numerator and denominator). Notice that
$(p/8)_4 = +1$ for primes $p \equiv 1 \bmod 16$ and  $(p/8)_4 = -1$ 
for primes $p \equiv 9 \bmod 16$. Moreover, $(\eps_1/p_2)$ is the 
quadratic residue character of $\eps_1 \bmod \fp$ (if 
$p_2 \equiv 1 \bmod 4$), where $\fp$ is a prime ideal in $k_1$ above 
$p_2$; for $\Delta_2 = 8$ and $\Delta_1 \equiv 1 \bmod 8$, the symbol 
$(\eps_1/8)$ is defined by
$(\eps_1/8) = (-1)^{T/4}$, where $\eps_1 = T+U\sqrt{\Delta_1}$.

\begin{cor} 
Let $p = \Delta_1$ and $q = \Delta_2 \equiv 1 \bmod 4$ be positive 
prime discriminants, and assume that $\Delta_2$ is fixed; then
$$ \begin{array}{lcccl}
 4 | h^+(k) & \iff & (\Delta_1/\Delta_2) = 1 & \iff  
                        & p \in \Spl(\Omega_4^+(\Delta_2)/\Q) \\
 4 | h(k) & \iff & (\Delta_1/\Delta_2)_4 = (\Delta_2/\Delta_1)_4 & \iff 
			& p \in \Spl(\Omega_4(\Delta_2)/\Q) \\
 8 | h^+(k) & \iff & (\Delta_1/\Delta_2)_4 = (\Delta_2/\Delta_1)_4 = 1 & \iff 
			& p \in \Spl(\Omega_4^+(\Delta_2)/\Q) \end{array} $$
\end{cor}
Here, the {\it governing fields} $\Omega_j(\Delta_2)$ are defined by
\begin{align*} 
 \Omega_4^+(\Delta_2) & = \ \Q(i,\sqrt{\Delta_2}\,), \\
 \Omega_4(\Delta_2)   & = \ \Omega_4^+(\sqrt{\eps_2}\,)
             \ = \ \Q(i,\sqrt{\Delta_2},\sqrt{\eps_2}\,), \\
 \Omega_8^+(\Delta_2) & = \Omega_4(\sqrt[4]{\Delta_2}\,) 
                \ = \ \Q(i,\root {4\,} \of {\Delta_2},\sqrt{\eps_2}\,). 
\end{align*}
The reason for studying governing fields comes from the fact 
that sets of primes splitting in a normal extension have 
Dirichlet densities.
The existence of fields governing the property $8|h^+(k)$ allows
us to conclude that there are infinitely many such fields. Governing
fields for the property $8|h(k)$ or $16|h^+(k)$ are not known and
conjectured not to exist. Nevertheless the primes $\Delta_1 = p$ such that
$8 \mid h(k)$ (or $16\mid h(k)$ etc.) appear to have exactly
the Dirichlet density one would expect if the corresponding
governing fields existed. 
Governing fields were introduced by Cohn and Lagarias \cite{CL81, CL83}
(see also Cohn's book \cite{Coh85}) and studied by
Morton \cite{Mor82a, Mor82b, Mor83, Mor90, Mor90a} and Stevenhagen
\cite{Ste88, Ste89, Ste93}. A typical result is

\begin{prop}
Let $p \equiv 1 \bmod 4$ and $r \equiv 3 \bmod 4$ be primes and
consider the quadratic number field $k = \Q(\sqrt{-rp}\,)$. Then
$8 \mid h(k) \iff (-r/p)_4 = +1$.
\end{prop}

More discussions on unramified cyclic quartic extensions of 
quadratic number fields can be found in Herz \cite{Her57}, 
Vaughan \cite{Vau85}, and Williams \& Liu \cite{WL94}. Another 
discussion of R\'edei's construction was given in Zink's 
dissertation \cite{Zin74}.

The class field theoretical approach was also used by R\'edei
\cite{Red43,Red53} (see Gerasim \cite{Geras}), as well as 
Furuta \cite{Fur}, Morton \cite{Mor90} and Benjamin, Lemmermeyer
\& Snyder \cite{BLS}. 

The problem for values divisible by squares was treated (along the 
lines of Dirichlet) by Perott \cite{Pero} and taken up again by 
R\'edei \cite{RedPl} and, via class field theory, by Jensen 
\cite{Jens1,Jens2,Jens3} and B\"ulow \cite{Bue}. Brown \cite{Bro} 
proved a very special case of Scholz's results using the theory of 
binary quadratic forms; see also Kaplan \cite{Kap1}. 
Buell \cite{Buell} gave a list of known criteria.

Despujols \cite{Des} showed that the norm of the fundamental unit 
is $(-1)^{t-r}$, where $t$ is the number of ramified primes, and
$r$ the number of ambiguous ideal classes containing an ambiguous 
ideal.

\section{Graphs of Quadratic Discriminants}

In this section we will explain the graph theoretical description
of classical results on the $4$-rank of class groups (in the strict
sense) of real quadratic number fields, and of solvability criteria
of the negative Pell equation. The connection with graph theory was
first described by Lagarias \cite{Lag} and used later by Cremona
\& Odoni \cite{CO}. Similar constructions were used by 
Vazzana \cite{Vaz1,Vaz2} for studying $K_2$ of the ring of integers 
$\cO_k$, as well as by Heath-Brown \cite{HB} and later by 
Feng \cite{Feng}, Li \& Tian \cite{LT}, and Zhao \cite{Zh1,Zh2} 
for describing the $2$-Selmer group of elliptic curves 
$Y^2 = X^3 - d^2X$.

\subsection*{The Language of Graphs}
A (nondirected) graph consists of a set $V$ of vertices and a subset 
$E \subseteq V \times V$ whose elements are called edges.

The degree of a vertex $d_i$ of a graph is the number of edges 
$(d_i,d_j) \in V$ adjacent to $d_i$. A graph is called Eulerian
if all vertices have even degree; graphs are Eulerian if and only
if there is a path through the graph passing each edge exactly once.

A tree is a connected graph containing no cycles (closed paths 
involving at least three vertices inside a graph). A subgraph 
of a graph $\gamma$ is called a spanning tree of $\gamma$ if it 
is a tree and if it contains all vertices.

\subsection*{Graphs of Quadratic Fields}
Let $d$ be the discriminant of a quadratic number field.
Then $d$ can be factored uniquely into prime discriminants:
$d = d_1 \cdots d_n$. Here the $d_i$ are discriminants of
quadratic fields in which only one prime $p_i$ is ramified.

Let $d$ be a discriminant of a real quadratic number field
such that all the $d_i$ are positive (equivalently, $d$ is
the sum of two integral squares). This implies, by quadratic
reciprocity, that $(d_i/p_j) = (d_j/p_i)$ for all $1 \le i, j \le n$.

To any discriminant $d$ as above we associate a graph 
$\gamma(d)$ as follows:
\begin{itemize}
\item $V = \{d_1, \ldots, d_n\}$;
\item $E = \{(d_i,d_j): (d_i/p_j) = (d_j/p_i) = -1\}$.
\end{itemize} 

Every factorization $d = \Delta_1 \Delta_2$ of $d$ into two
discriminants $\Delta_1, \Delta_2$ of quadratic fields
corresponds to a bipartitioning $\{A_1, A_2\}$ of the vertices  
by putting $A_1 = \{d_i: d_i \mid \Delta_1$ and
$A_2 = \{d_i: d_i \mid \Delta_2$. Clearly we have 
$A_1 \cup A_2 = V$ and $A_1 \cap A_2 = \varnothing$,
and if the factorization is nontrivial, we also have 
$A_1, A_2 \ne \varnothing$.

To each such bipartition, we associate a subgraph 
$\gamma(\Delta_1,\Delta_2)$ of $\gamma(d)$ by deleting all 
edges between vertices in $V_1$, and all edges between vertices 
in $V_2$; thus $\gamma(\Delta_1,\Delta_2)$ has vertices 
$V = V_1 \cup V_2$ and edges 
$E_{1,2} = \{(d_i,d_j) \in E: i \in V_1, j \in V_2\}$.

\begin{lem}
The factorization $d = \Delta_1\Delta_2$ is a $C_4$-decomposition
if and only if $\gamma(\Delta_1,\Delta_2)$ is Eulerian.
\end{lem}

\begin{proof}
The graph $\gamma(\Delta_1,\Delta_2)$ is Eulerian if and only
if for each $d_i \mid \Delta_1$ there is an even number of 
$d_j \mid \Delta_2$ such that $(d_j/p_i) = -1$. This is equivalent
to $(\Delta_2/p_i) = +1$. The claim follows.
\end{proof}

We call $\{A_1, A_2\}$ an Eulerian Vertex Decomposition (EVD) of
$\gamma(d)$ if the subgraph $\gamma(\Delta_1,\Delta_2)$ is Eulerian.
Since the number of $C_4$-decompositions equals the $4$-rank $e_2$ 
of the class group in the strict class group of $k = \Q(\sqrt{d}\,)$,
we see 

\begin{prop}
The number of EVDs of $\gamma(d)$ is $2^{e_2}$, where $e_2$ is the
$4$-rank of the class group of $k = \Q(\sqrt{d}\,)$ in the strict sense. 
\end{prop}

The R\'edei matrix can be interpreted as the adjacency matrix
of a graph $\Gamma(d)$; if $V = V_1 \cup V_2$ and 
$V_1 \cap V_2 = \varnothing$, the graph with the same vertices 
as $\Gamma(d)$ and the edges within $V_1$ and $V_2$ deleted
coincides with $\gamma(V_1,V_2)$.

The graph $\gamma(d)$ is said to be odd if it has the following
property: for every bipartitioning of $V$, there is an $a_1 \in A_1$
that is joined to an odd number of $a_2 \in A_2$, or vice versa.

\medskip

 [graph removed. maybe one day latex will be able
  to deal with jpg, gif, ps or pdf files]
\medskip

The graph $\gamma(5 \cdot 13 \cdot 17)$ is odd, the graph
$\gamma(5 \cdot 29 \cdot 41)$ is not. Observe that the 
negative Pell equation is solvable in the first, but not
in the second case.

\begin{prop}
Let the discriminant $d$ of a quadratic number field be a sum of
two squares. If $\gamma(d)$ is an odd graph, then $N \eps = -1$.
\end{prop}

\begin{proof}
Let $(x,y)$ be the solution of $X^2 - dY^2 = 1$ with minimal $y > 0$.
Then $x$ is odd, hence $x+1 = 2fr^2$, $x+1 = 2gs^2$, $fg = d$,
$1 = fr^2 - gs^2$. If $g = 1$, then $N\eps = -1$, and $f = 1$
contradicts the minimality of $y$. If $f, g > 1$ then we claim
that $1 = fr^2 - gs^2$ is not solvable in integers. Let 
$d = p_1 \cdots d_n$, $A = \{d_i: d_i \mid f\}$ and
$B= \{d_i: d_i \mid g\}$. Then $A \cup B = V = \{d_1, \ldots, d_n\}$,
$A \cap B = \varnothing$. Since $\gamma(d)$ is odd, we may
assume that there is an $d_i \in A$ that is adjacent to an odd 
number of $d_j \in B$. This implies $(g/p_i) = -1$, contradicting
the solvability of $1 = fr^2 - gs^2$.
\end{proof}

Lagarias observed that a congruence proved by Pumpl\"un \cite{Pum}
could be interpreted as follows:

\begin{prop}
Let $d$ be the discriminant of a quadratic number field. Then
$$ h^+(d) \equiv \sum_T \prod_{(d_i,d_j) \in T} 
            \Big(1 - \Big(\frac{d_i}{p_j}\Big)\Big) 
               \equiv 2^{n-1} \kappa_d \bmod 2^n, $$
where $T$ runs over all spanning trees of the complete graph 
with $n$ vertices, and where $\kappa_d$ is the number of 
spanning trees of $\gamma(d)$.
\end{prop}

This implies that $\Cl_2^+$ is elementary abelian if and only
if $\gamma(d)$ is an odd graph. This result is implicitly 
contained in Trotter \cite{Tro}.

\subsection*{Directed Graphs}
If $d = p_1 \cdots p_n$ is a product of primes $p_i \equiv 3 \bmod 4$,
then the graph with vertices $d_i = -p_i$ and adjacency matrix
$A = (a_{ij})$ defined by 
$$ \Big(\frac{p_j}{p_i}\Big) = \begin{cases}
         (-1)^{a_{ij}} & \text{if}\ i \ne j \\
         (-1)^{n+1} (\frac{d/p_i}{p_i}) & \text{if}\ i = j \end{cases} $$
is a directed graph (actually a tournament graph since each edge has 
a unique direction) studied by Kingan \cite{Kin}. From R\'edei's results,
Kingan deduced the following facts (see also Sueyoshi \cite{Sue}):
\begin{prop}
If $n$ is even, then $r_4(d) = n-1-\rank A$ and
$r_4(-d) = n-\rank A$ or $n-1-\rank A$.

If $n$ is odd, then $r_4(d) = n-1-\rank A$ or $r_4(d) = n-2-\rank A$, 
and $r_4(-d) = n-2-\rank A$.

Define $c_i$ via $(-1)^{c_i} = (2/p_i)$ and put $v = (c_1, \ldots, c_n)^T$.
Then in the cases where the rank formula is ambiguous, the greater value
is attained if $v \in \im (A-I)$.
\end{prop}

Kohno \& Nakahara \cite{KN} and Kohno, Kitamura, \& Nakahara 
\cite{KKN} used oriented graphs to describe Morton's results 
about governing fields and the computation of the $2$-part of 
the class group of quadratic fields.

Parts of the theory of R\'edei and Reichardt has been 
extended to function fields in one variable over finite fields; 
Ji \cite{Ji1} discussed decompositions of the second kind
over function fields, and in \cite{Ji2} he proved results of 
Trotter \cite{Tro} in this case, using the graph theoretic 
language discussed above.

\subsection*{Density Problems}
In this section we will address the question for how many $\Delta$
the negative Pell equation $X^2 - \Delta Y^2 = -1$ is solvable. It
was already noticed by Brahmagupta (see Whitford \cite{Whit}) that 
the solvability implies that $\Delta$ must be a sum of two squares.   

Consider therefore the set $\cD$ of quadratic discriminants not
divisible by any prime $\equiv 3 \bmod 4$, and let $\cD(-1)$ denote
the subset of all discriminants in $\cD$ for which the negative
Pell equation is solvable. The problem is then to determine whether
the limit
\begin{equation}\label{ESt}
   \lim_{x \to \infty} 
    \frac{\# \{\Delta \in \cD(-1): \Delta \le x\}}
         {\#    \{\Delta \in \cD : \Delta \le x\}} 
\end{equation}
exists, and if it does, to find it.

Such questions were first asked by Nagell \cite{Nag} and R\'edei 
\cite{RedMW,Redas}; using criteria for the solvability of the negative
Pell equation, R\'edei could prove that 
$$ \liminf_{x \to \infty} 
    \frac{\# \{\Delta \in \cD(-1): \Delta \le x\}}
         {\#    \{\Delta \in \cD : \Delta \le x\}} 
    > \alpha: = \prod_{j=1}^\infty (1 - 2^{1-2j}) = 0.419422\ldots,$$
a result later proved again by Cremona \& Odoni \cite{CO}. Stevenhagen
\cite{Ste93c} gave heuristic reasons why the density in (\ref{ESt})
should equal $1 - \alpha = 0.580577\ldots$.

Related questions concerning the density of quadratic fields 
whose class groups have given $4$-rank were studied by 
Gerth \cite{Ger1,Ger2} and Costa \& Gerth \cite{CG}.

Observe that the negative Pell equation is just one among
Legendre's equations; we might similarly ask for the density
of discriminants $d \equiv 1 \bmod 8$ for which $2x^2 - dy^2 = 1$
is solvable, among all discriminants $d \equiv 1 \bmod 8$ whose 
prime factors are $\equiv \pm 1 \bmod 8$. 

\subsection*{Selmer groups}
The graph theoretic language introduced for studying the density
of discriminants for which the negative Pell equation is solvable
was also employed for computing the size of Selmer groups of 
elliptic curves with a rational point of order $2$. See 
Heath-Brown \cite{HB}, Feng \cite{Feng}, Li \& Tian \cite{LT},
and Zhao \cite{Zh1,Zh2}.

\end{document}